\documentclass[a4paper,12pt]{amsart}

\usepackage{xcolor}
\usepackage{amsmath}
\usepackage{amssymb}
\usepackage{mathrsfs}
\usepackage{ifthen}
\usepackage{graphicx}
\usepackage[T1]{fontenc} 

\nonstopmode \numberwithin{equation}{section}
\newtheorem{thm}[equation]{Theorem}
\newtheorem{cor}[equation]{Corollary}
\newtheorem{lem}[equation]{Lemma}

\theoremstyle{definition}

\newtheorem{prob}[equation]{Problem}

\newcounter{minutes}
\divide\time by 60
\newcounter{hours}
\multiply\time by 60
\addtocounter{minutes}{-\time}

\newcounter {own}
\def\theown {\thesection       .\arabic{own}}

\newenvironment{pf}[1][]{%
 \vskip 3mm
 \noindent
 \ifthenelse{\equal{#1}{}}%
  {{\slshape Proof. }}%
  {{\slshape #1.} }%
 }%
{\qed\bigskip}

\newcounter{alphabet}

\def\be{\begin{equation}}
\def\ee{\end{equation}}

\newcommand{\bee}{\begin{enumerate}}
\newcommand{\eee}{\end{enumerate}}

\newcommand{\blem}{\begin{lem}}
\newcommand{\elem}{\end{lem}}
\newcommand{\bthm}{\begin{thm}}
\newcommand{\ethm}{\end{thm}}
\newcommand{\bcor}{\begin{cor}}
\newcommand{\ecor}{\end{cor}}
\newcommand{\beg}{\begin{examp}}
\newcommand{\eeg}{\end{examp}}
\newcommand{\begs}{\begin{examples}}
\newcommand{\eegs}{\end{examples}}
\newcommand{\bdefe}{\begin{defin}}
\newcommand{\edefe}{\end{defin}}
\newcommand{\bprob}{\begin{prob}}
\newcommand{\eprob}{\end{prob}}
\newcommand{\bei}{\begin{itemize}}
\newcommand{\eei}{\end{itemize}}

\begin{document}

\title{Support points of some classes of analytic and univalent functions}

\author{Vasudevarao Allu}
\address{Vasudevarao Allu,
School of Basic Sciences,
Indian Institute of Technology  Bhubaneswar,
Argul, Bhubaneswar, PIN-752050, Odisha (State),  India.}
\email{avrao@iitbbs.ac.in}

\author{Abhishek Pandey}
\address{Abhishek Pandey,
School of Basic Sciences,
Indian Institute of Technology  Bhubaneswar,
Argul, Bhubaneswar, PIN-752050, Odisha (State),  India.}
\email{ap57@iitbbs.ac.in}


\def\thefootnote{}
\footnotetext{ {\tiny File:~\jobname.tex,
printed: \number\year-\number\month-\number\day,
          \thehours.\ifnum\theminutes<10{0}\fi\theminutes }
} \makeatletter\def\thefootnote{\@arabic\c@footnote}\makeatother

\thanks{}

\maketitle
\pagestyle{myheadings}
\markboth{Vasudevarao  Allu and Abhishek Pandey}{Support points of some classes of analytic and univalent functions}

\begin{abstract}
Let $\mathcal{A}$ denote the class of analytic functions in the unit disk $\mathbb{D}:=\{z\in\mathbb{C}:|z|<1\}$ satisfying $f(0)=0$ and $f'(0)=1$.
Let $\mathcal{U}$ be the class of functions  $f\in\mathcal{A}$ satisfying 
$$\left|f'(z)\left(\frac{z}{f(z)}\right)^2-1 \right|< 1 \quad\mbox{ for } z\in\mathbb{D},$$
and $\mathscr{G}$ denote the class of functions $f\in \mathcal{A}$   satisfying $${\rm Re\,}\left(1+\frac{zf''(z)}{f'(z)}\right)>-\frac{1}{2} \quad\mbox{ for } z\in\mathbb{D}.$$
In the present paper, we characterize the set of support points of the classes $\mathcal{U}$ and $\mathscr{G}$. 


\end{abstract}

\section{Introduction and Preliminaries}\label{Introduction}

Let $\mathcal{H}$ denote the class of analytic functions in the unit disk $\mathbb{D}:=\{z\in\mathbb{C}:\, |z|<1\}$. Here  $\mathcal{H}$ is 
a locally convex topological vector space endowed with the topology of uniform convergence over compact subsets of $\mathbb{D}$. Let $\mathcal{A}$ denote the class of functions 
$f\in \mathcal{H}$ such that $f(0)=0$ and $f'(0)=1$. Let $\mathcal{S}$ 
denote the class of functions $f\in\mathcal{A}$ which are univalent ({\it i.e.}, one-to-one) in $\mathbb{D}$. 
Each function $f\in\mathcal{S}$ has the following representation
\begin{equation}\label{p2_eq1}
f(z)= z+\sum_{n=2}^{\infty}a_n z^n.
\end{equation}
A set $D \subset \mathbb{C}$ is said to be starlike with respect to a point $z_{0}\in D$ if the line segement joining $z_{0}$ to every other point $z\in D$ lies entirely in $D$. A set $D$ is called convex if line segment joining any two points of $D$ lies entirely in $D$. A function $f\in\mathcal{A}$ is called starlike (convex respectively) if $f(\mathbb{D})$ is starlike with respect to the origin (convex respectively). Denote by $\mathcal{S}^*$ and $\mathcal{C}$  the classes of starlike and convex functions in $\mathcal{S}$ respectively. It is well-known that a function $f\in\mathcal{A}$ belongs to $\mathcal{S}^*$ if,  and only if, ${\rm Re\,}\left(zf'(z)/f(z)\right)>0$ for $z\in\mathbb{D}$. Similarly, a function $f\in\mathcal{A}$ belongs to $\mathcal{C}$ if,  and only if, ${\rm Re\,}\left(1+zf''(z)/f'(z)\right)>0$ for $z\in\mathbb{D}$. Alexander \cite{Alexander-1915} proved that $f\in\mathcal{C}$ if, and only if, $zf'\in\mathcal{S}^*$. Given $\alpha\in(-\pi/2,\pi/2)$ and $g\in \mathcal{S}^*$, a function $f\in \mathcal{A}$ is said to be close-to-convex with argument $\alpha$  with respect to $g$ if 
$$ {\rm Re\,} \left(e^{i\alpha}\frac{zf'(z)}{g(z)}\right)>0 \mbox{ for } z\in \mathbb{D}.$$
Let $\mathcal{K}_{\alpha}(g)$ denote the class of all such functions, and $$\mathcal{K}(g):=\bigcup\limits_{\alpha\in(-\pi/2,\pi/2)}\,\mathcal{K}_{\alpha}(g) \mbox{ and } \mathcal{K}_{\alpha}:= \bigcup\limits_{g\in\mathcal{S}^*}\,\mathcal{K}_{\alpha}(g)$$
be the classes of close-to-convex functions with respect to $g$, and close-to-convex functions with argument $\alpha$, respectively. Let
$$\mathcal{K}:=\bigcup\limits_{\alpha\in(-\pi/2,\pi/2)}\,\mathcal{K}_{\alpha}=\bigcup\limits_{g\in \mathcal{S}^{*}}\, \mathcal{K}(g)$$
be the class of close-to-convex functions. It is well-known that every close-to-convex function is univalent in $\mathbb{D}$. Geometrically, $f\in \mathcal{K}$ means that the complement of the image domain $f(\mathbb{D})$ is the union of non-intersecting half- lines. These standard classes are related by the proper inclusions $\mathcal{C}\subsetneq\mathcal{S^*} \subsetneq\mathcal{K}\subsetneq\mathcal{S}$.\\

For $0<\lambda\le1$, let $\mathcal{U}(\lambda)$ be the class of  functions $f\in\mathcal{A}$ satisfying 
$$\left |f'(z)\left(\frac{z}{f(z)}\right)^{2}-1\right|<\lambda  \quad\mbox {  for }  z\in \mathbb{D}.
$$
Since $f'(z)(z/f(z))^{2}$ has only finite values, each function in the class $\mathcal{U}(\lambda)$ is non-vanishing in $\mathbb{D}\setminus\{0\}$.
Set $\mathcal{U}:=\mathcal{U}(1)$.  
It is also clear that the functions  in $\mathcal{U}(\lambda)$ are locally univalent. 
Furthermore, Aksentev \cite{Aksentev-1958} and Ozaki and Nunokawa \cite{ozaki-1972}  have shown that  functions in $ \mathcal{U}(\lambda)$ are univalent, {\it i.e.},
$\mathcal{U}(\lambda) \subseteq \mathcal{S}$  for  $0<\lambda\le 1$. Functions from the class $\mathcal{U}$ are univalent but not all are starlike which might be expected from the similarity in their analytic representations. This makes them interesting since the class of starlike functions is very large and in theory of univalent functions it is significant if a class does not entirely lie inside $\mathcal{S}^*.$ The Properties of meromorphic functions associated with the class $\mathcal{U}(\lambda)$ has been studied by Ali {\it et al.} in \cite{Ali-2020}. The integral mean problem and arc length problem for functions in $\mathcal{U}$ have been studied in \cite{Ali-2020}. As we know that neither $\mathcal{U}\subset\mathcal{S}^*$ nor $\mathcal{S}^* \subset \mathcal{U}$ but Allu and Pandey \cite{Allu-2020} have proved that $\overline{co}\,\mathcal{U}=\overline{co}\,\mathcal{S}^*$. The main observation was that $z/(1-xz)^{2}\in \mathcal{S}^*\cap\,\mathcal{U}$
for each $x$ such that $|x|=1$. For details of the class $\mathcal{U}(\lambda)$ we refer to \cite[Chapter 12]{Vasu-book} \\

Let $\mathscr{G}$ denote the class of functions $f\in \mathcal{A}$   satisfying $${\rm Re\,}\left(1+\frac{zf''(z)}{f'(z)}\right)>-\frac{1}{2} \quad\mbox{ for } z\in\mathbb{D}.$$
Functions in $\mathscr{G}$ are convex in one direction (see \cite{Umezama-1952}) and hence are close-to-convex  but they are not necessarily starlike in $\mathbb{D}$. The importance of class $\mathscr{G}$ in the case of certain univalent harmonic mappings is exposed in \cite{Bshouty-1986}. The properties of sections of functions in $\mathscr{G}$  have been studied in \cite{Bharanedhar-2014} and \cite{Ponnusamy-2014}. Moreover, if $f\in \mathscr{G}$ and is of the form (\ref{p2_eq1}), then we have the following coefficient inequality
\begin{equation}\label{p2_eq2} 
|a_{n}|\le \frac{n+1}{2} \,\,\,\mbox{ for all } n\ge 2
\end{equation}

and the equality in \ref{p2_eq2} holds only for the function $g_{0}(z)$ and its rotations, where
$$g_{0}(z)=\frac{z-(\frac{x}{2})z^{2}}{(1-xz)^{2}},\quad |x|=1.$$\\
Suppose $X$ is a linear topological vector space and $V \subseteq X$.  A point $x\in V$ is called an extreme point of $V$ if it has no representation of the form $x=ty+(1-t)z, 0<t<1$ 
as a proper convex combination of two distinct points $y$, $z\in V$. We denote by $EV$ the set of extreme points of $V$. The convex hull of a set $V \subseteq X$ is the smallest convex set containing $V$. The closed convex 
hull, denoted by ${\overline{co}}V$, is defined as the intersection of all closed convex sets containing $V$.  
Therefore  the closed convex hull of $V$ is the smallest closed convex set containing $V,$ which  is the closure of the convex hull of $V$. The Krein-Milman theorem \cite{Dunford-1958} asserts that every compact subset of a locally convex topological vector space is contained in the closed convex hull of its extreme points.\\
A function $f$ is called support point of a compact subset $\mathcal{F}$ of $\mathcal{H}$ if $f\in \mathcal{F}$ and if there is a continuous linear functional $J$ on $\mathcal{H}$ such that ${\rm Re\,}\,J$ is non-constant on $\mathcal{F}$ and
$${\rm Re\,}J(f)=\max\{{\rm Re\,}J(g):g\in\mathcal{F}\},$$
set of all support points of a compact family $\mathcal{F}$ is denoted by $supp\,\mathcal{F}$\\

 Support points of the families of starlike, convex functions have been studied in \cite{MacGregor-1971} and support points of close-to-convex functions have been studied in \cite{Grassmann-1976} and \cite{Wilken-1984}. Support points for starlike and convex function of order $\alpha$ have been studied in \cite{Cochrane-1978}.  
For a general reference, and for many important results on the above  topic, we refer to \cite{Hallenbeck-MacGregor-1984}. For a very recent work related to support points and extreme points we refer to \cite{Ponnusamy-2020}. In this paper Deng, Ponnusamy and Qiao has been proved a necessary and sufficient condition for harmonic Bloch mapping $f$ to be a support point of the unit ball of the normalized harmonic Bloch spaces in $\mathbb{D}$.\\

 Before going to prove our main results we mention some important lemmas which play a vital role to prove our main results.

\begin{lem}\cite[Theorem 4.3]{Hallenbeck-MacGregor-1984}\label{lem1}
$J$ is a complex-valued continuous linear functional on $\mathcal{H}$ if, and only if, there is a sequence $\{b_{n}\}$ of complex numbers satisfying
$$\overline{\lim\limits_{n\to\infty}}|b_{n}|^{1/n}<1$$
and such that
$$J(f)=\sum_{n=0}^{\infty}b_{n}a_{n}$$
where $f\in \mathcal{H}$ and $f(z)=\sum_{n=0}^{\infty}a_{n}z^{n}$ for $|z|<1$.
\end{lem}

\begin{lem}{\cite[Theorem 4.5]{Hallenbeck-MacGregor-1984}}\label{lem2}
Let $\mathcal{F}$ be a compact subset of $\mathcal{H}$ and $J$ be a complex-valued continuous linear functional on $\mathcal{H}$. Then
$\max\{ {\rm Re\, }J(f): f\in \overline{co}\, \mathcal{F}\} = \max\{ {\rm Re\, }J(f): f\in \mathcal{F}\} = \max\{ {\rm Re\, }J(f): f\in E\overline{co}\, \mathcal{F}\}$.
\end{lem}

\begin{lem}\cite{Allu-2020}
$\overline{co}\, \mathcal{U}$ consists of all functions represented by  
\begin{equation*}
f(z)=\int_{|x|=1}\frac{z}{(1-xz)^{2}}\,\, d\mu(x), 
\end{equation*}
where $\mu \in \wedge$, and $\wedge$ denotes the set of probability measures on $\partial{\mathbb{D}}$.
Further,  $E\overline{co}\, \mathcal{U}$ consists of the functions of the form 
\begin{equation*}
f(z)=\frac{z}{(1-xz)^{2}},\quad |x|=1.
\end{equation*}
\end{lem}

\begin{lem}\cite{Muhanna-2014}
$\overline{co}\, \mathscr{G}$ consists of all functions represented by  
\begin{equation*}
f(z)=\int_{|x|=1}\frac{z-(\frac{x}{2})z^{2}}{(1-xz)^{2}}\,\, d\mu(x), 
\end{equation*}
where $\mu \in \wedge$, and $\wedge$ denotes the set of probability measures on $\partial{\mathbb{D}}$.
Further,  $E\overline{co}\, \mathscr{G}$ consists of the functions of the form 
\begin{equation}\label{p2-eq2}
f(z)=\frac{z-(\frac{x}{2})z^{2}}{(1-xz)^{2}},\quad |x|=1.
\end{equation}
\end{lem}
 
The main aim of this paper is to characterize the set of support points of the class $\mathcal{U}$ and $\mathscr{G}$.

\section{Support points of the class $\mathcal{U}$}
\begin{thm}\label{supp}
$supp\,\mathcal{U}=E\overline{co}\,\mathcal{U}=\left\{\frac{z}{(1-xz)^{2}}:|x|=1\right\}.$
\end{thm}
\begin{pf}
First we show that the functions of the form $z/(1-xz)^{2}$ where $|x|=1$, belong to $supp\,\mathcal{U}$. For $x\in \mathbb{C}$ such that $|x|=1$, consider the function
$$f_{0}(z)=\frac{z}{(1-xz)^{2}}=z+\sum_{n=2}^{\infty}\,nx^{n-1}z^{n}=z+2xz^{2}+\sum_{n=3}^{\infty}nx^{n-1}z^{n}.$$
A simple computation shows that $f_{0}\in \mathcal{U}$.
Now, consider the continuous linear functional $\phi$ defined by, 
$\phi(f)=\bar{x}a_{2}$. Clearly, ${\rm Re\, }\phi$ is non-constant on $\mathcal{U}$. We know that $\phi(f_{0})=\bar{x}2x=2$ which gives that ${\rm Re\, }\phi(f_{0})=2$. That is
$${\rm Re\, }\phi(f_{0})=\max\{{\rm Re\, } \phi(f):f\in\mathcal{U}\}.$$
Hence each function of the form $z/(1-xz)^{2}$ is a support point of $\mathcal{U}$. That is
\begin{equation}\label{eq1}
E\overline{co}\,\mathcal{U}\subseteq supp\, \mathcal{U}.
\end{equation}
Our next aim is to show that $supp\,\mathcal{U}\subseteq E\overline{co}\,\mathcal{U}$. In order to prove this, we first prove that each function in the class $supp\,\overline{co}\,\mathcal{U}$ is of the form 
\begin{equation}\label{eq}
\sum_{k=1}^{m}\lambda_{k}\frac{z}{(1-x_{k}z)^{2}}
\end{equation} where $\lambda_{k}\ge 0$, $\sum_{k=1}^{m}\lambda_{k}=1$ and $|x_{k}|=1$ and $m=1,2,3,\ldots.$\\

Let $\widetilde{f}$ be any support point of $\overline{co}\,\mathcal{U}$. That is there exists a continuous linear functional on $\mathcal{H}$ say $\widetilde{J}$ with ${\rm Re\, }\widetilde{J}$ is non-constant on $\overline{co}\,\mathcal{U}$ such that 
$${\rm Re\, }\widetilde{J}(\widetilde{f})=\max\{{\rm Re\, }\widetilde{J}(f):f\in\overline{co}\,\mathcal{U}\}.$$
By Lemma \ref{lem1}, we obtain a sequence $\{b_{n}\}$ such that $\limsup\limits_{n\to\infty}|b_{n}|^{1/n}<1$ 
and
$$\widetilde{J}(f)=\sum_{n=0}^{\infty}b_{n}a_{n}$$
where $f\in \mathcal{H}$ and $f(z)=\sum_{n=0}^{\infty}a_{n}z^{n}$ for $|z|<1$. Consider 
$$F(z,x)=\frac{z}{(1-xz)^{2}}=z+\sum_{n=2}^{\infty}a_{n}(x)z^{n}$$
where $a_{n}(x)=nx^{n-1}$. Then $|a_{n}(x)|\le nr^{n-1}$ if $1/r<|x|<r$ for $r>1$.\\
Let ${\rm Re\, }\widetilde{J}(\widetilde{f})=\widetilde{M}$. Hence in view of Lemma \ref{lem2}, we obtain
\begin{eqnarray}
&{\rm Re\, }\widetilde{J}(\tilde{f})=&
\max\{{\rm Re\, }\widetilde{J}(f):f\in\overline{co}\,\mathcal{U}\} =\max\{{\rm Re\, }\widetilde{J}(f):f\in\mathcal{U}\}\\ \nonumber
&&:=\max\{{\rm Re\, }\widetilde{J}(f):f\in E\overline{co}\,\mathcal{U}\}=\widetilde{M}.
\end{eqnarray}
Let $G(x)=\widetilde{J}(F(z,x))=b_{1}+\sum_{n=2}^{\infty}b_{n}a_{n}(x).$ Since $\limsup\limits_{n\to\infty}|b_{n}|^{1/n}<1$ for some $\alpha$ satisfying $0<\alpha<1$ and some $N\in\mathbb{N}$, we obtain 
$|b_{n}|^{1/n}\le\alpha$ whenever $n>N$. It follows that $|b_{n}a_{n}(x)|\le \alpha^{n}nr^{n-1}=(n/r)(\alpha r)^{n}$. Choosing $r$ such that $\alpha r<1$, one can see that the series representation of $G(x)$ is convergent if $1/r<|x|<r$. Thus, $G$ is analytic function on $\{x:1/r<|x|<r\}$. Let $H$ be defined by 
$$H(x)=\frac{1}{2}\left(G(x)+\overline{G\left(\frac{1}{\bar{x}}\right)}\right).$$
It is easy to observe that $H$ is analytic on $\{x:1/r<|x|<r\}$. Also we note that if $|x|=1$ then 
\begin{eqnarray*}
& H(x)=& \frac{1}{2}\left(G(x)+\overline{G\left(\frac{1}{\bar{x}}\right)}\right)=\frac{1}{2}\left(G(x)+\overline{G\left(\frac{x}{|x|^{2}}\right)}\right)\\ \nonumber &&=\frac{1}{2}\left(G(x)+\overline{G(x)}\right)={\rm Re\, }G(x).
\end{eqnarray*}

Suppose $H(x)=\widetilde{M}$ has infinitely many solutions on $\{x:|x|=1\}$. Since $\{x:|x|=1\}$ is a compact set and $H$ is analytic on $\{x:1/r<|x|<r\}$ we obtain $$H(x)=\widetilde{M}\mbox{ for } x \mbox{ on } \{x:1/r<|x|<r\}.$$
In particular, $$H(x)=\widetilde{M}\mbox{ for all } x\in \{x:|x|=1\},$$
That is
$${\rm Re\, }G(x)=\widetilde{M} \mbox{ for all } x\in \{x:|x|=1\}$$
which gives
$${\rm Re\, }\widetilde{J}(F(z,x))=\widetilde{M} \mbox{ for all } x\in \{x:|x|=1\}.$$
Now ${\rm Re\, }\widetilde{J}$ is constant on $E\overline{co}\,\mathcal{U}$ is evident from the fact that $E\overline{co}\,\mathcal{U}=\{F(z,x):|x|=1\}$. It follows that ${\rm Re\, }\widetilde{J}$ is constant on $\overline{co}\,\mathcal{U}$ which is not possible. Therefore $H(x)=\widetilde{M}$ has only finitely many solutions on $\{x:|x|=1\}$. Hence,
$$H(x)=\widetilde{M} \mbox{ has only finitely many solutions on } \{x:|x|=1\}.$$
That is
$$ {\rm Re\, }\widetilde{J}(F(z,x))=\widetilde{M} \mbox{ for finitely many points on } \{x:|x|=1\}.$$
That is, there are only finitely many functions in $E\overline{co}\,\mathcal{U}$ which maximizes ${\rm Re\, }\widetilde{J}$ over $\overline{co}\,\mathcal{U}$.\\

Now consider the set of all functions in $\overline{co}\,\mathcal{U}$ which maximizes ${\rm Re\, }\widetilde{J}$ over $\overline{co}\,\mathcal{U}$. That is, consider the following set $$\mathscr{C}=\{g:g\in\overline{co}\,\mathcal{U} \mbox{ and } {\rm Re\, }\widetilde{J}(g)=\max\{{\rm Re\, }\widetilde{J}(f):f\in\overline{co}\,\mathcal{U}\}\}.$$
It is easy to see that $\mathscr{C}$ is convex and compact. Hence by Krein-Milman theorem $E\mathscr{C}$ is non-empty. Also suppose there is a function $f^*\in E\mathscr{C}$ but $f^*\notin E\overline{co}\,\mathcal{U}$, then there are $f_{1}^{*}$ and $f_{2}^{*}$ in $\overline{co}\,\mathcal{U}$ such that
$$f^*(z)=\beta f_{1}^{*}(z)+(1-\beta)f_{2}^{*}(z) \mbox{ for some } \beta\in (0,1).$$ 
By the linearity of ${\rm Re\, }\widetilde{J}$, we obtain
$${\rm Re\, }\widetilde{J}(f^*(z))=\beta\,{\rm Re\, }\widetilde{J}(f_{1}^{*}(z))+(1-\beta)\,{\rm Re\, }\widetilde{J}(f_{2}^{*}(z)).$$
This shows that $f_{1}^{*}$ and $f_{2}^{*}$ belong to $\mathscr{C}$ which is a contradiction since $f^*\in E\mathscr{C}$. Thus, $E\mathscr{C}\subset E\overline{co}\,\mathcal{U}$. Since $E\overline{co}\,\mathcal{U}$ contains only finite number of elements which are also in $E\mathscr{C}$, it follows that $E\mathscr{C}$ is a finite set of distinct functions $f_{k}$ of the form 
$$f_{k}(z)=\frac{z}{(1-x_{k}z)^{2}}, \quad\mbox{ where } |x_{k}|=1.$$
Hence 
$$\mathscr{C}=\left\{f:f(z)=\sum_{k=1}^{m}\lambda_{k}\frac{z}{(1-x_{k}z)^{2}},\,\, \lambda_{k}\ge0,\,\, \sum_{k=1}^{m}\lambda_{k}=1 \mbox{ and } |x_{k}|=1 \right\}.$$
From the definition of $\mathscr{C}$, it is clear that $\widetilde{f}\in \mathscr{C}$ and so $\widetilde{f}$ is of the form (\ref{eq}). Thus we have proved that every support point of $\overline{co}\,\mathcal{U}$ is of the form (\ref{eq}).\\

Our next aim is to show that $supp\,\mathcal{U}\subseteq E\overline{co}\,\mathcal{U}$.
For this, let $\widehat{f}\in\mathcal{U}$ be any support point of $\mathcal{U}$. Since $\mathcal{U}\subset \overline{co}\,\mathcal{U}$ therefore $\widehat{f}\in \overline{co}\,\mathcal{U}$, there is continuous linear functional $\widehat{J}$ on $\mathcal{H}$ with ${\rm Re\, }\widehat{J}$ is a non-constant on $\mathcal{U}$ such that 
$${\rm Re\, }\widehat{J}(\widehat{f})=\max\{{\rm Re\, }\widehat{J}(f):f\in\mathcal{U} \}.$$
Since $\mathcal{U}$ is compact, in view of Lemma \ref{lem2} we obtain
$${\rm Re\, }\widehat{J}(\widehat{f})=\max\{{\rm Re\, }\widehat{J}(f):f\in\overline{co}\,\mathcal{U} \}.$$
Therefore $\widehat{f}$ is a support point of $\overline{co}\,\mathcal{U}$. Hence it is of the form (\ref{eq}). That is,
$$\widehat{f}=\sum_{k=1}^{m}\lambda_{k}\frac{z}{(1-x_{k}z)^{2}}.$$
If $\lambda_{k}\ne 0$ for at least two values of $k$ then $\widehat{f}$ has at least two poles at $\bar{x_{k}}$ on the unit circle each of order $2$ and such function can not be univalent on $\mathbb{D}$. Hence, 
$$\widehat{f}(z)=\frac{z}{(1-xz)^{2}}\,\,\,\, \mbox{ where } |x|=1.$$
Thus, $\widehat{f}\in E\overline{co}\,\mathcal{U}$ which gives
\begin{equation}\label{eq3}
supp\,\mathcal{U}\subseteq E\overline{co}\,\mathcal{U}
\end{equation}

and hence by (\ref{eq1}) and (\ref{eq3}), we obtain
$$supp\,\mathcal{U}=E\overline{co}\,\mathcal{U}=\left\{\frac{z}{(1-xz)^{2}}:|x|=1\right\}.$$
This completes the proof.
\end{pf}
\section{Support point of the class $\mathscr{G}$}
\begin{thm}
$supp\,\mathscr{G}=E\overline{co}\,\mathscr{G}=\left\{\displaystyle\frac{z-\frac{x}{2}z^{2}}{(1-xz)^{2}}:|x|=1\right\}.$
\end{thm}
\begin{pf}
First we show that each function of the form (\ref{p2-eq2}) is a support point of $\mathscr{G}$. For $x\in \mathbb{C}$ such that $|x|=1$, consider the function
$$g_{0}(z)=\frac{z-(\frac{x}{2})z^{2}}{(1-xz)^{2}}=z+\sum_{n=2}^{\infty}\,\frac{n+1}{2}x^{n-1}z^{n}=z+\frac{3}{2}xz^{2}+\sum_{n=3}^{\infty}\frac{n+1}{2}x^{n-1}z^{n}.$$ A simple computation shows that for each $x$ such that $|x|=1$, $g_{0}\in \mathscr{G}$. Consider the continuous linear functional $\phi:\mathcal{H}\to \mathbb{C}$ defined by $\phi(f)=\bar{x}a_{2}$. Then $\phi(g_{0})=3/2$ which gives that ${\rm Re\, }\phi(g_{0})=3/2$. Therefore from (\ref{p2_eq2}) we obtain,
$${\rm Re\, }\phi(g_{0})=\max\{{\rm Re\, } \phi(g):g\in\mathscr{G}\}.$$
Therefore each function of the form $(z-(x/2)z^{2})/(1-xz)^{2}$ is a support point of $\mathscr{G}$. That is
\begin{equation}\label{eq4}
E\overline{co}\,\mathscr{G}\subseteq supp\, \mathscr{G}.
\end{equation}
Our next aim is to show that $supp\,\mathscr{G}\subseteq E\overline{co}\,\mathscr{G}$. In order to prove this, we first prove that each function in the class $supp\,\overline{co}\,\mathscr{G}$ is of the form 
\begin{equation}\label{eq5}
\sum_{k=1}^{m}\lambda_{k}\frac{z-\frac{x_{k}}{2}z^{2}}{(1-x_{k}z)^{2}}
\end{equation} where $\lambda_{k}\ge 0$, $\sum_{k=1}^{m}\lambda_{k}=1$ and $|x_{k}|=1$ for $m=1,2,3,\ldots$.\\

 Let $\widetilde{g}$ be any support point of $\overline{co}\,\mathscr{G}$. That is there exists a continuous linear functional on $\mathcal{H}$ say $\widetilde{L}$ with ${\rm Re\, }\widetilde{L}$ is non-constant on $\overline{co}\,\mathscr{G}$  such that 
$${\rm Re\, }\widetilde{L}(\widetilde{g})=\max\{{\rm Re\, }\widetilde{L}(g):g\in\overline{co}\,\mathcal{U}\}.$$
By Lemma \ref{lem1}, we obtain a sequence $\{c_{n}\}$ such that $\limsup\limits_{n\to\infty}|c_{n}|^{1/n}<1$
and
$$\widetilde{L}(f)=\sum_{n=0}^{\infty}b_{n}a_{n}$$
where $f\in \mathcal{H}$ and $f(z)=\sum_{n=0}^{\infty}a_{n}z^{n}$ for $|z|<1$. Consider 
$$\widetilde{F}(z,x)=\frac{z-\frac{x}{2}z^{2}}{(1-xz)^{2}}=z+\sum_{n=2}^{\infty}a_{n}(x)z^{n}$$
where $$a_{n}(x)=\frac{n+1}{2}x^{n-1},$$ which implies that, $$|a_{n}(x)|\le \frac{n+1}{2}r^{n-1} \mbox{ for }\,\, 1/r<|x|<r,\, r>1.$$
Let ${\rm Re\, }\widetilde{L}(\widetilde{g})=\widehat{M}$. Then in view of Lemma \ref{lem2}, we obtain
\begin{eqnarray}
&{\rm Re\, }\widetilde{L}(\widetilde{g})=&
\max\{{\rm Re\, }\widetilde{L}(g):g\in\overline{co}\,\mathscr{G}\} =\max\{{\rm Re\, }\widetilde{L}(g):g\in\mathscr{G}\}\\ \nonumber
&&=\max\{{\rm Re\, }\widetilde{L}(g):g\in E\overline{co}\,\mathscr{G}\}=\widehat{M}.
\end{eqnarray}
Let  $$\widetilde{G}(x)=\widetilde{L}(\widetilde{F}(z,x))=c_{1}+\sum_{n=2}^{\infty}c_{n}\frac{n+1}{2}x^{n-1}(x).$$ Since $\limsup\limits_{n\to\infty}|c_{n}|^{1/n}<1$, we can easily show that the series representation of $\widetilde{G}(x)$ is convergent if $1/r<|x|<r$ for $r$ sufficiently close to $1$. Thus, $\widetilde{G}$ is an analytic function on $\{x:1/r<|x|<r\}$. Let $\widetilde{H}$ be defined by  
$$\widetilde{H}(x)=\frac{1}{2}\left(\widetilde{G}(x)+\overline{\widetilde{G}\left(\frac{1}{\bar{x}}\right)}\right).$$
Clearly, $\widetilde{H}$ is an analytic function on $\{x:1/r<|x|<r\}$. We note that if $|x|=1$ then 
$\widetilde{H}(x)={\rm Re\, }\widetilde{G}(x).$ By using the similar argument in the proof of Theorem \ref{supp}, we obtain $$\widetilde{H}(x)=\widehat{M} \mbox{ has finitely many solutions on } \{x:|x|=1\}.$$ That is
$$ {\rm Re\, }\widetilde{L}(\widetilde{F}(z,x))=\widehat{M} \mbox{ for finitely many points in } \{x:|x|=1\}.$$
That is, there are only finitely many functions in $E\overline{co}\,\mathscr{G}$ which maximizes ${\rm Re\, }\widetilde{L}$ over $\overline{co}\,\mathscr{G}$.\\

Consider the set of all functions in $\overline{co}\,\mathscr{G}$ which maximizes ${\rm Re\, }\widetilde{L}$ over $\overline{co}\mathscr{G}$. That is consider the set
$$\mathscr{K}=\{g:g\in\overline{co}\,\mathscr{G}\mbox{ and } {\rm Re\, }\widetilde{L}(g)=\max\{{\rm Re\, }\widetilde{L}(f):f\in\overline{co}\,\mathscr{G}\}\}.$$
It is easy to see that $\mathscr{K}$ is convex and compact. Therefore by Krein-Milman theorem $E\mathscr{K}$ is non-empty. Also one can observe that $E\mathscr{K}\subset E\overline{co}\,\mathscr{G}$. It follows that since there are only finite number of elements in $E\overline{co}\,\mathscr{G}$ which are also in $E\mathscr{K}$, it follows that $E\mathscr{K}$ is a finite set of distinct functions $g_{k}$ of the form 
$$g_{k}(z)=\frac{z-\frac{x_{k}}{2}z^{2}}{(1-x_{k}z)^{2}} \quad\mbox{ where } |x_{k}|=1.$$
Hence 
$$\mathscr{K}=\left\{f:f(z)=\sum_{k=1}^{m}\lambda_{k}\frac{z-\frac{x}{2}z^{2}}{(1-x_{k}z)^{2}},\,\, \lambda_{k}\ge0,\,\, \sum_{k=1}^{m}\lambda_{k}=1 \mbox{ and } |x_{k}|=1 \right\}.$$
From the definition of $\mathscr{K}$, it is clear that $\widetilde{g}\in \mathscr{K}$ and hence $\widetilde{g}$ is of the form (\ref{eq5}). Thus we have proved that every support point of $\overline{co}\,\mathscr{G}$ is of the form (\ref{eq5}).\\

Finally, we show that $supp\,\mathscr{G}\subseteq E\overline{co}\,\mathscr{G}$.
For this, let $\widehat{g}\in\mathscr{G}$ be any support point of $\mathscr{G}$. Since $\mathscr{G}\subset \overline{co}\,\mathscr{G}$, we have $\widehat{g}\in \overline{co}\,\mathscr{G}$. Therefore, there exists a continuous linear functional $\widehat{\mathcal{L}}$ on $\mathcal{H}$ with ${\rm Re\, }\widehat{\mathcal{L}}$ is a non-constant on $\mathscr{G}$ such that 
$${\rm Re\, }\widehat{\mathcal{L}}(\widehat{g})=\max\{{\rm Re\, }\widehat{\mathcal{L}}(g):g\in\mathscr{G} \}.$$
Since $\mathscr{G}$ is compact, in view of Lemma \ref{lem2}, we obtain
$${\rm Re\, }\widehat{\mathcal{L}}(\widehat{g})=\max\{{\rm Re\, }\widehat{\mathcal{L}}(g):g\in\overline{co}\,\mathscr{G} \}.$$
It can be seen that $\widehat{g}$ is a support point of $\overline{co}\,\mathscr{G}$. Therefore it is of the form (\ref{eq5}). That is
$$\sum_{k=1}^{m}\lambda_{k}\frac{z-(\frac{x_{k}}{2})z^{2}}{(1-x_{k}z)^{2}}.$$
If $\lambda_{k}\ne 0$ for at least two values of $k$ then $\hat{g}$ has at least two poles at $\bar{x_{k}}$ on the unit circle each of order $2$ and such function can not be univalent on $\mathbb{D}$. Therefore, 
$$\widehat{g}(z)=\frac{z-(\frac{x}{2})z^{2}}{(1-xz)^{2}}\,\,\,\, \mbox{ where } |x|=1.$$
Thus, $\widehat{g}\in E\overline{co}\,\mathscr{G}$, which shows that 
\begin{equation}\label{eq6}
supp\,\mathscr{G}\subseteq E\overline{co}\,\mathscr{G}.
\end{equation}
Hence by (\ref{eq4}) and (\ref{eq6}), we obtain
$$supp\,\mathscr{G}=E\overline{co}\,\mathscr{G}=\left\{\frac{z-\frac{x}{2}z^{2}}{(1-xz)^{2}}:|x|=1\right\}.$$
This completes the proof.
\end{pf}

\vspace{4mm}
\noindent\textbf{Acknowledgement:} 
The first author thank SERB-MATRICS 
and the second author thank PMRF-MHRD, Govt. of India for their support.


\begin{thebibliography}{99}
    

\bibitem{Muhanna-2014}
{\sc Y. Abu-Muhanna}, {\sc Liulan Li} and {\sc S. Ponnusamy}, Extremal problems on the class of convex functions of order $-1/2$, {\it Arch. Math.} {\bf 103} (2014), 461--471.    
    
    
\bibitem{Aksentev-1958}
{\sc L.A. Akesentev}, Sufficient conditions for univalence of certain integral representations (Russian), {\it Izv. Vyss. Ucebn. Zavrd. Matematika} {\bf 4} (1958), 3--7.

 \bibitem{Alexander-1915}
 {\sc Alexander J. W.}, Functions which maps the interior of the unit circle upon simple regions. {\it Ann. of Math.}, {\bf 17} (1915--1916), 12--22.
 
\bibitem{Ali-2020}
{\sc Md Firoz Ali}, {\sc A. Vasudevarao} and {\sc  Hiroshi Yanagihara}, On a class of univalent functions defined by a differential inequality, {\it J. Ramanujan Math. Soc}  (To appear)	arXiv:2006.15577 (2020).
 

 
\bibitem{Allu-2020}
{\sc A. Vasudevarao} and {\sc A. Pandey}, The Zalcman conjecture for certain analytic and univalent functions, {\it J. Math. Anal. Appl.} {\bf 492} (2), 2020. 
  
\bibitem{Bharanedhar-2014}
{\sc S. V. Bharanedhar} and {\sc S. Ponnusamy}, Uniform close-to-convexity radius of sections of functions in the close-to-convex family, {\it J. Ramanujan Math. Soc.} {\bf 29} (2014), 243--251. 
  
  \bibitem{MacGregor-1971}
   {\sc L. Brickman} and {\sc T.H. MacGregor} and {\sc D. R. Wilken} Convex hull of some classical family of univalent functions, {\it Trans. Amer. Math. Soc.} {\bf 156} (1971), 91--107.  
  
\bibitem{Bshouty-1986}
{\sc D. Bshouty} and {\sc A. Lyzzaik}, close-to-convexity criteria for planar harmonic mappings, {\it Complex Anal. Oper. Theory}, {\bf 5} (2011), 767--774.
    
\bibitem{Cochrane-1978}
    {\sc  P. C. Cochrane} and {\sc T. H. MacGregor}, {$\rm Fr\acute{e}chet$} differentiable functionals and suppprt points for families of analytic functions,{\it Trans. Amer. Math. Soc.} {\bf 236}, (1978), 75--92.   
    
\bibitem{Ponnusamy-2020}
{\sc Hua Deng}, {\sc Saminathan Ponnusamy} and {\sc Jinjing Qiao}, Extreme Points and Support Points of Families of Harmonic Bloch Mappings, {\it Potential Anal}(2020) https://doi.org/10.1007/s11118-020-09871-3.
    
\bibitem{Dunford-1958}
{\sc  N. Dunford} and {\sc J. T. Schwartz}, Linear operators, part I. Interscience: {\it New York}, 1958
    
\bibitem{Grassmann-1976}  
{\sc E. Grassman}, {\sc W. Hengartner} and {\sc G. Schober}, Support points of the class of close-to-convex functions, {\it Canad. Math. Bull.} {\bf 19} (2), 1976    

 

\bibitem{Hallenbeck-MacGregor-1984}
{\sc D. J. Hallenbeck} and {\sc T. H. MacGregor}, Linear problem and convexity techniques in geometric function theory, Pitman, 1984.

\bibitem{Hallenbeck-1989}
 {\sc  D.J. Hallenbeck}, {\sc S. Perera} and {\sc D.R.Wilken}, Subordination, Extreme points and support points, {\it Complex Variables}, {\bf 11}, (1989), 111--124.

 
 
\bibitem{ozaki-1972}
{\sc S. Ozaki} and {\sc M. Nunokawa}, The Schwarzian derivative and univalent functions, {\it Proc. Amer. Math. Soc.} {\bf 33} (1972), 392--394.

\bibitem{Ponnusamy-2014}
{\sc S. Ponnusamy}, {\sc S. K. Sahoo} and {\sc H. Yanagihara}, Radius of convexity of partial sums of functions in the close-to-convex family, {\it Nonlinear Anal.} {\bf 95} (2014), 219--228. 
 
 \bibitem{Vasu-book}
 {\sc Derek K. Thomas}, {\sc Nikola Tuneski} and {\sc Allu Vasudevarao}, Univalent functions. A primer, De Gruyter Studies in Mathematics, {\bf 69}. De Gruyter, Berlin, 2018.
 
\bibitem{Umezama-1952}
{\sc T. Umezama}, Analytic functions convex in one direction, {\it J. Math. Soc. Japan}, {\bf 4} (1952), 194--202.

\bibitem{Wilken-1984}
{\sc D. Wilken}, {\sc R. Hornblower}, On the support points of close-to-convex functions, {\it Houston Journal Of Mathematics}, {\bf 10} (4), 1984.
\end{thebibliography}
\end{document}